\documentclass[24pt,reqno]{article}
\usepackage
[colorlinks=true, linkcolor=webgreen, filecolor=webbrown,
citecolor=webgreen]{hyperref}

\def\R{{\mathbb R}}

\def\C{{\mathbb C}}
\usepackage{amsmath}
\usepackage{amssymb}
\usepackage{color}
\usepackage{epsf}
\usepackage{graphics}

\begin{document}
\makeatletter

\begin{center}
\epsfxsize=10in
\end{center}

\def\endofproofmark{$\Box$}

\begin{center}
\vskip 1cm {\LARGE\bf On Periodic Solutions of Li\'{e}nard
Equation} \vskip 1cm \vspace{10mm}

{\large Ali Taghavi}

\vskip .5cm
Department of Mathematical Sciences\\
Sharif University of Technology\\
P.O. Box 11365­9415, Tehran, Iran \vskip 3mm
School of Mathematics\\
Institute for Studies in Theoretical Physics and Mathematics\\
P.O. Box 19395­1795, Tehran, Iran\vskip 3mm {\tt taghavi@ipm.ir}
\end{center}

\newtheorem{theorem}{Theorem}
\newtheorem{prop}{Proposition}
\newtheorem{lemma}{Lemma}
\newtheorem{cor}{Corollary}
\newtheorem{prob}{Note and Problem}
\newtheorem{conj}{Conjecture}
\newtheorem{ques}{Question}
\def\frameqed{\framebox(5.2,6.2){}}
\def\deshqed{\dashbox{2.71}(3.5,9.00){}}
\def\ruleqed{\rule{5.25\unitlength}{9.75\unitlength}}
\def\myqed{\rule{8.00\unitlength}{12.00\unitlength}}
\def\qed{\hbox{\hskip 6pt\vrule width 7pt height11pt depth1pt\hskip 3pt}
\bigskip}
\newenvironment{proof}{\trivlist\item[\hskip\labelsep{\bf Proof}:]}{\hfill
 $\frameqed$ \endtrivlist}
\newcommand{\COM}[2]{{#1\choose#2}}

\thispagestyle{empty} \null \addtolength{\textheight}{1cm}

\begin{abstract}
It is conjectured by Pugh, Lins and de Melo in [7] that the system
of equations
$$
\begin{cases}
\dot{x}=y-F(x)\\
\dot{y}=-x
\end{cases}
$$
has at most $n$ limit cycles when the degree of $F=2n+1$ or $2n+
2$. Put $M$ for uniform upper bound of the number of limit cycles
of all systems of equations of the form
$$
\begin{cases}
\dot{x}=y -(ax^4 + bx^3 + cx^2 + dx)\\
\dot{y}=-x.
\end{cases}
$$
In this article, we show that $M\neq 2$. In fact, if an example
with two limit cycle existed, one could give not only an example
with $n+2$ limit cycles for the first system, but also one could
give a counterexample to the conjecture $N(2,3)=2$ [see the
conjecture $N(2,3)=2$ of F. Dumortier and C.Li, Quadratic
Li\'{e}nard equation with Quadratic Damping, J. Differential
Equation, 139 (1997) 41­59]. We will also pose a question about
complete integrability of Hamiltonian systems in $\R^4$ which
naturally arise from planner Li\'{e}nard equation. Finally,
considering the Li\'{e}nard equation as a complex differential
equation, we suggest a related problem which is a particular case
of conjecture. We also observe that the Li\'{e}nard vector fields
have often trivial centralizers among polynomial vector fields.
\end{abstract}
\section{Introduction}
We consider Li\'{e}nard equation in the form
\begin{eqnarray}
\begin{cases}
\dot{x}=y-F(x)\\
\dot{y}=-x
\end{cases}
\end{eqnarray}
where $F(x)$ is a polynomial of degree $2n+1$ or $2n+2$. It is
conjectured in [7] that this system has at most $n$ limit cycles.
In particular, it was conjectured that the system
\begin{eqnarray}
\begin{cases}
\dot{x}=y-(ax^4+bx^3+cx^2+dx)\\
\dot{y}=-x
\end{cases}
\end{eqnarray}
has at most one limit cycle. The phase portrait of (1) is
presented in [7]. (1) has a center at the origin if and only if
$F(x)$ is an even polynomial. The following useful Lemma is proved
in [7]:
\begin{lemma}
Let $F(x)=E(x)+O(x)$ where $E$ is an even polynomial and $O$ is an
odd polynomial, and that $O(x)=0$ has a unique root at $x=0$. Then
system (1) does not have a closed orbit.
\end{lemma}
To prove the Lemma, it was shown that a first integral of the
system
\begin{eqnarray}
\begin{cases}
\dot{x}=y-E(x)\\
\dot{y}=-x
\end{cases}
\end{eqnarray}
that is analytic and defined on $\R^2-{0}$ is a monotone function
along the solutions of (1). We introduce the following conjecture
about the above first integral:\\
\begin{conj} Let $E(x)$ be an even polynomial of degree at
least 4, then there is no global analytic first integral for
system (3) on $\R^2$.
\end{conj}
The reason for conjecture: There are two candidates for defining a
first integral, namely the square of the intersection of the
solution with negative $y$­axis and the other the square of the
intersection of the solution with $x$­axis. The first is well
defined on $\R^2$ but certainly is not analytic at the origin and
the second is analytic in the region of closed orbits but it can
not be defined on all of $\R^2$. In fact, there are solutions not
intersecting the $x$­axis: Consider the region surrounded by $y =
\frac{-1}{x}$ and $y=0$ for $x\ggg1$ and look at the direction
field of
$$
\begin{cases}
\dot{x}=y-F(x)\\
\dot{y}=-x
\end{cases}
$$
on the boundary of this region. We conclude that there is a
solution remaining in this area for all $t<0$ whenever the
solution is defined. The analyticity of the second function,
intersection with $x$­axis, is obtained by the fact that the
1­from $(y - E(x))dy+xdx$ is the pull back of $(y-g(x))dy+dx$
under $\pi(x,y)=(x^2,y)$ with $E(x)=g(x^2)$. Now, the differential
equation
$$
\begin{cases}
\dot{x}=y-g(x)\\
\dot{y}=-1
\end{cases}
$$
does not have a singular point and the intersection of the orbits
with $x$­axis defines an analytic function as a first integral,
which we call $K(x,y)$. Then $K(x^2,y)$ is a first integral for
the original system (3). Note that putting $E(x)=x^4$, then the
dual from of (3) is the pull back of Ricati equation, whose
solution can not be determined in terms of elementary functions
(This can be seen using Galois Theory, see [6]). I thank professor
R. Roussarie for hinting me to this pull back. In line of the
above conjecture, we propose the following question:
\begin{ques} Let an analytic vector field on the plane have a
non­degenerate center. As a rule, is it possible to define
analytic first integrals globally in the domain of periodic
solutions? (However, there is one locally.) In fact, using Riemann
mapping theorem we can assume that the region of the center is all
of $\R^2$.
\end{ques}
We return to Li\'{e}nard equation (1). The limit cycles of (1)
correspond to fixed points of Poincar\'{e} return map. Let $F(x)$
has odd degree with positive leading coefficient. As a rule, we
can not define Poincar\'{e} return map on positive $y$­semi­axis.
For example, if $F(x)=x^{2n+1}+2x$, in which case the singular
point is a node, there is no solution starting on positive
$y$­axis and returning again to this axis; see the direction field
$$
\begin{cases}
\dot{x}=y-(x^{2n+1}+2x)\\
\dot{y}=-x
\end{cases}
$$
on the semi­line $y=x$ and $x>0$. Then, contrary to what is
written in [7] or is pointed out in [10], we can define a
Poincar\'{e} return map only in the case that origin is a weak or
strong focus or in the case of existence of at least one limit
cycle. In [10], it is also pointed out that Dulac's problem is
trivial for (1), that is, for any given $F(x)$, (1) has a finite
number of limit cycles. Let $F(x)$ has odd degree, then
$P_0=\lim_{y\rightarrow +\infty}P(y)$ exists and even we have
$\lim_{y\rightarrow +\infty}(P(y)-P_0)y^i=0$ for all $i$. Then $P$
has a finite number of fixed points. (In this work we do not use
the strong approach $\lim_{y\rightarrow +\infty}y^i(P(y)-P_0)=0$
and thus we do not prove it.) For the case that $F$ has even
degree, the above is not trivial and one can deduce it from the
results of this paper. In fact we must consider the case that
$\lim_{y\rightarrow +\infty}P(y)=+\infty$. Equivalently, we have a
loop at Poincar\'{e} sphere based at infinity. As a simple
consequence of the above Lemma, we note that the system
$$
\begin{cases}
\dot{x}=y-x^2\\
\dot{y}=\epsilon(a-x)
\end{cases}
$$
does not has limit cycles for $a\neq 0$, because putting $x:=x+a$,
$y:=y-a^2$ we will obtain
$$
\begin{cases}
\dot{x}=y-(x^2+2ax)\\
\dot{y}=-\epsilon x.
\end{cases}
$$
Using the Lemma, this system does not have a limit cycle. In
figures on page 478 and 479 of [4], it appears that the existence
of limit cycles is claimed.
\section{Main Results}
\begin{theorem} For any $a,~b,~c$, there exists a unique $d=(a,b,c)$ such
that the system (2) has a homoclinic loop in Poincar\'{e} sphere.
This loop is stable if and only if the singular point is unstable.
\end{theorem}
\begin{cor} The maximum number of limit cycles of (2) can not
be exactly two.
\end{cor}
\begin{proof} (Proof of the theorem) For fixed $a,~b,~c,$ $a>0$, $b>0$,
put $U(d)$ and $S(d)$ for intersection of unstable and stable
manifolds corresponding to the topological saddle $(0,1,0)$ on the
equator of Poincar\'{e} sphere. In fact $U(d)$ and $S(d)$ are
similar to $P_{-}$ and $P_{+}$ on page 339 of [7], figure 3.
$U(d)$ and $S(d)$ are continuous and monotone functions (as will
be proved). We also prove $U(d)=S(d)$ has a unique root. First
note that for $d>0$ there is no limit cycle or homoclinic loop
based at infinity, using the lemma. Then we assume $d<0$, note
that $U(0)>S(0)$, otherwise from the stability of the origin in
(2) for $b>0$ and $d=0$, we would have a limit cycle for (2) which
is a contradiction with the lemma. On the other hand, $U(d)\lll
S(d)$ for $d\lll -1$ because the minimum value of $F(x)=ax^4+bx^3+
cx^2+dx$ goes to $-\infty$ as $d\rightarrow -\infty$ and $S$ is
decreasing. Now by continuity and monotonicity of $U$ and $S$ we
have a unique root $d_0=\psi(a,b,c, )$ such that $U(d_0)=S(d_0)$.
We must prove that $U$ and $S$ are continuous and monotone. It can
be directly observed, even without using any classical theorem,
that $U$, $S$ are continuous. We prove that $S$ is decreasing,
similarly $U$ is increasing. Let $\gamma_d$ be a solution of
$(2)_d$ that is asymptotic to the graph of
$F(x)=ax^4+bx^3+cx^2+dx$, in $x<0$ in fact, $\gamma_d$ is the
stable manifold for the topological saddle $(0,1,0)$ that
intersects $y$­axis in $S(d)$. We prove that for $d'<d<0$,
$S(d')>S(d)$. $\gamma_d$ is a curve without contact for $(2)_{d'}$
and the direction field of $(2)_{d'}$ on $\gamma_d$ is toward
''left'' and the direction of $(2)_{d'}$ on the semi­line $y>0$.
$x=0$, is toward ''right''. Then there is a unique orbit of
$(2)_{d'}$ that remains in the region surrounded by $\gamma_d$ and
the semi­line $x=0$, $y>0$. This orbit is a stable manifold for
$(0,1,0)$ of system $(2)_{d'}$, say; $\gamma_{d'}$ certainly
$\gamma'_{d}$ can not intersect $\gamma_d$ so $\gamma_{d'}$ will
intersect negative $y$­axis in $S(d')$ above $S(d)$. Therefore $S$
is decreasing. For the proof of the theorem it remains to prove
that the loop is attractive, still assuming $a>0,~b>0,~d<0$.
Before proving that the homoclinic loop is attractive we point out
that, at first glance, this loop has a degenerate vertex at
$(0,1,0)$, i.e. the linear part of vector field at the vertex is $\left[%
\begin{array}{cc}
  0 & 0 \\
  0 & 0 \\
\end{array}%
\right]$. But using weighted compactification as explained in [3],
one has an elementary polycycle with two vertices, one vertex is a
''From'' vertex and another is a ``To''. Thus we do not have an
``unbalanced polycycle'' and therefore behavior of solution near
the polycycle can not be easily determined. Thus, we use the
following direct computation. (for definition of balanced and
unbalanced polycycles, see [5] page 21.) Let
$P_0=U(d_0)=S(d_0)<0$; We have a Poincar\'{e} return map $P$
defined on negative vertical section $(P_0,0)$, we will prove that
$\lim_{y\rightarrow P_0} P'(y)=0$, then $P(y)-y$ is negative for
$y$ near $P_0$. Note that considering the orbit of points near
$P_0$ in the vertical section $(P_0,0)$ is equivalent to
considering the orbits of points $(0,\tilde{y})$ with
$\tilde{y}\ggg 1$. Let $\gamma$ be the orbit of (2) corresponding
to a homoclinic loop whose existence is proved above. $\gamma$ is
asymptotic to the graph of $F(x)=ax^4+bx^3+cx^2+dx$.

\begin{figure}[h]
\begin{center}
\resizebox{2.5in}{!}{\includegraphics{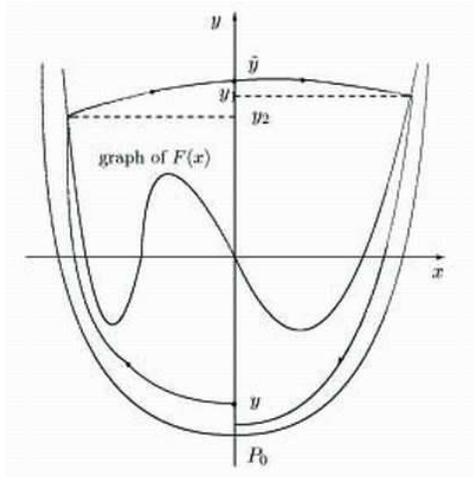}}
\end{center}
\caption{$F(x)=ax^4+bx^3+cx^2+dx$}
\end{figure}


The orbit starting from $(0,\tilde{y})$ intersects the graph of
$F$ at points with $y$­coordinates $y_1$ and $y_2$, when traced in
positive and negative time direction, respectively. Recall that we
want to prove $\lim_{y\rightarrow P_0} P'(y)=0$; we use the
following three facts:\\
\textbf{I)} There is a constant $K$ such that
$|\tilde{y}-y_1|<K\sqrt{\tilde{y}}$ and
$|\tilde{y}-y_2|<K\sqrt{\tilde{y}}$, where $K$ depends only on $a,
b, c, d$.\\
\textbf{II)} $\gamma$ is asymptotic to the graph of $F$.\\
\textbf{III)} Let $A(y)$ and $B(y)$ two right inverse of
$F(x)=ax^4+bx^3+cx^2+dx$, that is $F(A(y))=F(B(y))=y$ for $y\ggg
1$. Then $\lim_{y\rightarrow +\infty}(A(y)+B(y))=\frac{-b}{2a}$.
(III) is a simple exercise (II) is pointed out in [7]: So we only
prove the first. It suffices to prove (I), with the same notation
$y_1$ and $y_2$ for the intersection of the orbit with the graph
of $F(x)+1$ (in place of
$F(x)$).$|\tilde{y}-y_1|=\frac{|\tilde{y}-y_1|}{\Delta x}\Delta
x$, where $\Delta x$ is the $x$­coordinate of the intersection of
the orbit starting from $(0,\tilde{y})$ with the graph of
$F(x)+1$, then $\Delta x<K\sqrt{\tilde{y}}$ for some $K$, since
the degree of $F$ is 4 and $\frac{\tilde{y}-y_1}{\Delta
x}=\frac{-x}{y_0-f(x_0)}$, where $(x_0,y_0)$ is a point of the
trajectory starting at $(0,\tilde{y})$ and $0<x<\Delta x$ (using
mean value theorem). Then $\frac{\tilde{y}-y_1}{\Delta
x}<K\sqrt{\tilde{y}}$ and (I) is proved. Now we compute
$\lim_{y\rightarrow P_0} P'(y)$ (or equivalently for
$\tilde{y}\ggg 1$ as stated above): $P'(y)=\frac{y}{P(y)}
exp[\int_0^{T(y)}``div''dt]$ where by ``div'' we mean the
divergence of vector field (2). i.e. $-F'(x(t))$, and $T(y)$ is
the time of first return of the solution starting $(0,y)$ to
negative vertical section (for the standard formula of $P'(y)$ see
[1]).
$$
\int div dt=-\int
F'(x)dt=\int\frac{F'(x)}{x}=\int\frac{-F'(x)dr}{y-F(x)}
$$
We use $\tilde{y}$ for the intersection of the solution starting
at $(0,y)$ with the positive $y$­axis in positive time. We show
that $\int_0^{T(y)}div=\int_0^{T(y)}-F'(x(t))dt$ goes to $-\infty$
for $y$ near $P_0$. We divide $\int_0^{T(y)}divdt$ into three
parts $I_1$, $I_2$, $I_3$: $I_1$ for the part of the integral that
the orbit $\gamma_y$, the solution starting at $(0,y)$ for
$y\in(p_0,0)$ lice in $|x|<c$, where $c$ is a nonzero constant.
$I_2$ corresponds to the part of $\gamma_y$ above the horizontal
line $y=\tilde{y}-K\sqrt{\tilde{y}}$ where $K$ is the same
constant which is given in (I), and $I_3$ is the remaining part of
$I=\int_0^{T(y)}div$. In fact, in $I_3$ we compute the integral of
the divergence of (2) along the part of $\gamma_y$ that lies below
the horizontal line $y=\tilde{y}-K\sqrt{\tilde{y}}$ and outside of
$|x|<c_1$ where in this part $\gamma_y$ lies between the graph of
$F(x)$ and the orbit $\gamma$ which corresponds to the homoclinic
loop. Note that as $\tilde{y}\rightarrow +\infty$ and consequently
$y\rightarrow P_0$, $I_2$ and $I_3$ are very large in absolute
value but $I_1$ remains bounded because $\int
F'(x)=\int\frac{F'(x)dx}{y-F(x)}$. On the other hand we will see
that $I_2+I_3$ goes to $-\infty$ as $y\rightarrow P_0$, thus we
may ignore the term $I_1$. We could choose $K$ in (I) so that
$|\frac{F'(A(\tilde{y}))}{A(\tilde{y})}|<K\sqrt{\tilde{y}}$ where
$A(\tilde{y})$ is the positive inverse of $F(x)$ for large
$\tilde{y}$. There is a constant $K_1$ such that
$|I_2|<K_1\tilde{y}$. Now we show $\frac{I_3}{\tilde{y}}$ to
$-\infty$ as $y\rightarrow P_0$. Recall that $\tilde{y}$ is as in
the figure. Applying III, we realize that $I_3$ is the integral of
some function which is big in absolute value as $y\rightarrow P_0$
(or $\tilde{y}\rightarrow +\infty$). Generally speaking, let
$g(\gamma)$ be a function such that $\lim_{y\rightarrow +\infty}
g(\gamma)=+\infty$. Put $G(\gamma)=\int_0^{\gamma}g(s)ds$, then
$\lim_{y\rightarrow +\infty}\frac{G(\gamma)}{\gamma}=+\infty$. So,
it suffices to prove that $I_3$ is the integral of some function
with respect to $\tilde{y}$, where this function goes to $-\infty$
as $\tilde{y}\rightarrow +\infty$.
\begin{eqnarray}
I_3=\int\frac{F'(x)}{x}dy,
\end{eqnarray}
$$
\frac{F'(x)}{x}=4ax^2+3bx+2c+dx.
$$
We consider two parts of $\gamma_y$, one in $x>c$ and another in
$x<-c_1$ since below the horizontal line
$y=\tilde{y}-K\sqrt{\tilde{y}}$ the orbit $\gamma_y$ lies between
$\gamma$ and graph of $F(x)$. Note that $\gamma$ is asymptotic to
the graph of $F(x)$. Now, applying III, we will obtain
$\frac{I_3}{y}\rightarrow -\infty$: for large values of
$\tilde{y}_1$ the term $``\frac{d}{x}''$ in (1) can be omitted,
then we must compute $-\lim_{\tilde{y}+\infty}
4a(A^2(\tilde{y})-B^2(\tilde{y}))+3b(A(\tilde{y})-B(\tilde{y}))$,
where $A(\tilde{y})$ and $B(\tilde{y})$ are inverses of $F(x)$
such that $F(A(\tilde{y}))=F(B(\tilde{y}))=\tilde{y}$ and
$A(\tilde{y})>0$, $B(\tilde{y})<0$. We look at
$-\lim_{\tilde{y}+\infty}
(A(\tilde{y})-B(\tilde{y}))[4a(A(\tilde{y})-B(\tilde{y}))+3b]$.
Certainly $(A(\tilde{y})-B(\tilde{y}))\rightarrow +\infty$. So,
the above limit goes to $-\infty$ and
$\frac{I_3}{\tilde{y}}\rightarrow +\infty$ as $y\rightarrow P_0$.
This completes the proof of the theorem.
\end{proof}
\begin{proof} (Proof of the corollary) From the proof of the theorem, we
conclude that (2) can have an even number of limit cycles for
$d\leq d_0=\psi(a,b,c)$ and can have odd number of limit cycles
for $d_0<d<0$. Now let $(2)_d$ have exactly two limit cycles.
therefore $d\leq d_0=\psi(a,b,c)$. If $d<d_0$, then $(2)_{d_0}$
has at least two limit cycles, counting multiplicity. It is
because any closed orbit of $(2)_d$ is a curve without contact for
$(2)_{d_0}$ and the direction fields of $(2)_{d_0}$ on this closed
curve are toward the interior. But $(2)_{d_0}$ can not have an odd
number of limit cycles. For instance let none of the two limit
cycles of $(2)_{d_0}$ be semi­stable. Then by a small perturbation
$d=d_0-\epsilon$, (2); these two limit cycles do not die. On the
other hand, we would obtain another limit cycle near the loop
$\gamma$. Because for $d=d_0-\epsilon$, $\gamma$ is a curve
without contact for $(2)_d$ and the direction of $(2)_d$ on
$\gamma$ is toward the interior. Note that $\gamma$ was an
attractive loop, therefore if $\epsilon$ is very small, using
Poincar\'{e}­Bendixson Theorem we will obtain a third limit cycle
near the loop. Now any semi­stable limit cycle can be replaced by
two limit cycles by an appropriate perturbation $d=d_0±\epsilon$
depending on which side of the limit cycle is stable. For instant,
let $(2)_{d_0}$ have a semi­stable limit cycle in interior of the
loop then $d=d_0-\epsilon$ gives us two limit cycle near the
semi­stable ones and simultaneously one limit cycle near the loop.
However the corollary is proved, we point out that the existence
of 3 limit cycles for (2) easily implies the existence of 3 limit
cycles for
$$
\begin{cases}
\dot{x}=y-(ax^4+bx^3+cx^2+dx)\\
\dot{y}=-x+\epsilon x^2
\end{cases}
$$
for small $\epsilon$. Now by a linear change of coordinates we put
$\epsilon=1$, so we would obtain counterexample to the conjecture
that the latter system for $\epsilon=1$ has at most 2 limit
cycles. See conjecture $N(2,3)=2$ in [3].
\end{proof}
\textbf{Remark 1.} For the proof of the existence of the loop
$\gamma$ and also the proof of the corollary, in fact we used the
rotational property of the parameter $d$ in (2), namely, any
solution of $(2)_d$ is a curve without contact for $(2)_{d'}$, if
$d'\neq d$. In particular, periodic solutions of $(2)_d$ are
closed curves without contact for $(2)_{d'}$. The simple but
useful phenomenon of ``rotated vector field theory'' introduced by
Duff, is some­ times used erroneously. See, for example , the
investigation of
$$
\begin{cases}
\dot{x}=y-(ax^4+bx^3+dx)\\
\dot{y}=-x
\end{cases}
$$
in [8]. It is claimed that there is no limit cycle for
$d<d_0=\psi(a,b,c)<0$, where $b>0$. In fact, the following
situation that could occur, is not considered in [8]:\\
As we assumed above, let $a>0,~b>0$, for small $d<0$ we have
exactly one (small) Hopf bifurcating limit cycle. It is possible
that this limit cycle, before arriving to a loop situation, dies
out in a semi­stable limit cycle. Put $X={d|(5)_d {\rm has exactly
one limit cycle}}$, we do not necessarily have $d_0=infX$, where
$d_0=\psi(a,b,0)$, correspond to loop situation. Put ``i'' for the
above infimum. It could be $i>d_0$ and $(5)_i$ possesses a
semi­stable limit cycle. It is also possible that when the
outermost limit cycle is dying out in the loop, the two innermost
limit cycles have not arrived to each other yet. In fact [8]
suggests an affirmative answer to the conjecture for system
(2).\\\\
\textbf{Remark 2.} The homoclinic loop $\gamma$, as an orbit on
the plane, not on the Poincar\'{e} sphere, divides the plane into
two parts: its interior, where all solutions are complete, and its
exterior, where all solutions have a finite interval of
definition. Interior orbits are complete because $\gamma$ is a
complete orbit by virtue of $\int dt=\int\frac{dx}{y-f(x)}$ and
$\gamma$ being asymptotic to the graph of $F(x)$. The exterior
points of $\gamma$ are not complete orbits because they tend to
hyperbolic sink and the source on the equator of the Poincar\'{e}
sphere (See [2]). Now, a trivial observation is that Li\'{e}nard
equation (1) can not have an isochronous center i.e. a center with
a fixed period for all closed orbits surrounding it.\\\\
\textbf{Remark 3.} As we saw in the proof of the theorem and the
corollary, the inequality $U(d)<S(d)$ or the reverse, determines
oddness or evenness of the number of limit cycle. In this
direction we point out that: Let $F(x)$ be an even degree
polynomial with positive leading coefficient, and $U(\epsilon)$
and $S(\epsilon)$ be similar to $U(d)$ and $S(d)$ above for the
system
$$
\begin{cases}
\dot{x}=y-F(x)\\
\dot{y}=-\epsilon x
\end{cases}
$$
Then $\lim_{\epsilon\rightarrow 0}U(\epsilon)=M^+$ and
$\lim_{\epsilon\rightarrow 0}S(\epsilon)=M^-$ where $M^-$ and
$M^+$ are minimum values of $F(x)$ on $(-\infty,0]$ and
$[0,+\infty)$ resp. The proof is identical to the proof of
existence of orbits passing through $U(d)$ and $S(d)$ asymptotic
to the graph of $F(x)$, as in [7].\\\\
\textbf{Remark 4.} Note that when the degree of $F(x)$ in (1) is
odd, the behavior of infinity is determined only by the sign of
the leading term of $F(x)$. Then, giving an example of (2) with 3
limit cycles would give, inductively, $n+2$ limit cycles for (1),
for all n.\\\\
\textbf{Remark 5.} Considering ``flow'' version of the problem of
``centralizer of diffeomorphisms'' described in [10] one can
easily observe the following partial result. Let $L$ be the
Li\'{e}nard vector field similar to (1) with at least one closed
orbit and $X$ be a polynomial vector field such that
$[L,X]\equiv0$, then $X=cL$ where $c$ is a constant real number.
In general, let two vector fields have commuting flows and
$\gamma$ be a closed orbit for one of them which does not lie in
an isochronous band of closed orbits. Then $\gamma$ must be
invariant by another vector field and if both vector fields are
polynomials. Then either $\gamma$ is an algebraic curve or two
vector fields are constant multiple of each other. But Li\'{e}nard
systems do not have algebraic solutions, [9]. More generally by
the following proposition we have ``Non existence of algebraic
solution implies triviality of centralizer''.
\begin{prop} Let $M$ be the set of all polynomial vector fields on $\C^2$.
Then $X\in M$ has trivial centralizer if $X$ dose not have an
algebraic solution.
\end{prop}
\begin{proof} Let $[X,Y]=0$. Then $X.Det(X,Y)=(DivX).(Det(X,Y))$, so
$``Det(X,Y)=0''$ defines an algebraic curve invariant under $X$
(By $Det(X,Y)$) we mean the determinant of a 2×2 matrix whose
columns are components of $X$ and $Y$).
\end{proof}
\textbf{Remark 6.} The following could be a (real) generalization
of formula used in proof of proposition:\\
Let $(M,\omega)$ be a $2n$ dimensional symplectic manifold and $X,
Y$ two vector fields with the condition $[X,Y]\equiv 0$, then
$X.\omega(X,Y)=n(DivX)\omega(X,Y)$. This formula is trivial for
usual symplectic structure of $\R^2$. Further there is a local
chart around each point of a two dimensional symplectic manifold
that $\omega$ can be represented in the trivial form. Thus the
formula is proved for arbitrary two dimensional symplectic
manifold. Now, in general case we have a two dimensional
symplectic submanifold $N$ of $M$, that $X$ and $Y$ are tangent to
$N$. (For points $m\in M$ that $\omega(X(m),Y(m))\neq 0$, using
Frobenius theorem.) From other hand $DivX_M=nDivX_N$. The
investigation of points $m$ that $\omega(X(m),Y(m))=0$ is trivial.
Then the proof is completed.
\begin{ques} Dirac introduced the following embedding of planner vector fields into Hamiltonian
system in $\R^4: H=zP(x,y)+wQ(x, y)$. Now we consider the
Hamiltonian $H=z(y-F(x))-wx$. When is this Hamiltonian completely
integrable?
\end{ques}
By completely­integrable Hamiltonian, we mean that there is a
first integral for the system
$$
\begin{cases}
\dot{x}=y-F(x)\\
\dot{y}=-x\\
\dot{z}=w-F'(x)z\\
\dot{w}=-z
\end{cases}
$$
independent of $H$.\\
The particular case $F(x)=x^2$ for which (1) has a center with
global first integral $\phi(x,y)=(y-x^2+\frac{1}{2})e^{-2y}$,
suggests that when (1) does not have a center, the above
Hamiltonian is not completely integrable. In fact, the function
$\phi$, as above, is a first integral, independent of $H$ for
above four dimensional system. Then, when (1) is not integrable,
one expects that the corresponding Hamiltonian is not completely
integrable. However, surprisingly, putting $F(x)=kx$, we have
another first integral $yw+xz$.

\begin{ques} Considering (1) as a vector field on $\C^2$, and in line of
conjecture in [7] one can think of the validity of the following
two statement.\\
I)There are at most n different leaves containing real limit
cycles.\\
II)There are at most n real limit cycles lying on the same leaf.
By ``leaf'' we mean a leaf of the foliation corresponding to the
equation (1) on $\C^2-{0}$.
\end{ques}

\textbf{Acknowledgement.} The author thanks Professor S.
Shahshahani for fruitful conversations. He also acknowledges the
financial support of the Institute for Studies in Theoretical
Physics and Mathematics.

\end{document}